\documentclass[12pt]{amsart}
\newtheorem{theorem}{Theorem}[section]
\newtheorem{lemma}[theorem]{Lemma}
\newtheorem{proposition}[theorem]{Proposition}

\newtheorem{example}[theorem]{Example}

\newtheorem{corollary}[theorem]{Corollary}

\newcommand{\Co}{\mbox{$\mathbb{C}$}}
\newcommand{\N}{\mbox{$\mathbb{N}$}}

\newcommand{\E}{\mbox{${\mathcal E}$}}

\newcommand{\R}{\mbox{${\mathcal R}$}}

\newcommand{\Li}{\mbox{${\mathcal L}$}}
\newcommand{\M}{\mbox{${\mathcal M}$}}

\newcommand{\T}{\mbox{${\mathcal T}$}}

\begin{document}
   \title[The ideal envelope of an operator algebra]
{The ideal envelope of an operator algebra}
\author{David P. Blecher and Masayoshi Kaneda} 

\thanks{*This research was 
supported by a grant from the National Science Foundation}  
\thanks{Revision of Oct 8, 2001. }  
\address{Department of Mathematics\\University of Houston\\
4800 Calhoun\\Houston, TX 77204-3008, U.S.A.} 
\email{dblecher@math.uh.edu and kaneda@math.uh.edu}\maketitle

\vspace{10 mm}
 
\begin{abstract}
A left ideal of any $C^*$-algebra is an example of an 
operator algebra with a right contractive approximate 
identity (r.c.a.i.).   Conversely, we show here
and in \cite{Bid1} that 
operator algebras with r.c.a.i. should be studied in
terms of a certain left ideal of a $C^*$-algebra.   
We study operator algebras and 
their multiplier algebras from the perspective of `Hamana
theory' and using the multiplier algebras introduced by 
the first author.  
\end{abstract}

\pagebreak
\newpage

\section{Introduction and notation}
A left ideal of any $C^*$-algebra is an example of an
operator algebra with a right contractive approximate
identity.
 Conversely, we study
operator algebras with  right contractive approximate identity in
terms of a certain left ideal of a $C^*$-algebra. 
 
A (concrete) operator algebra is a subalgebra of $B(H)$,
for some Hilbert space $H$.   More abstractly,
an operator algebra will be an algebra $A$ with
a norm defined on the space $M_n(A)$ of $n \times n$ 
matrices with entries in $A$, for each $n \in \N$, such that
there exists a {\em completely isometric}\footnote{A map 
$T : X \rightarrow Y$ is 
completely isometric if $[x_{ij}] \mapsto
[T(x_{ij})]$ is isometric on $M_n(X)$ for all $n \in \N$} homomorphism
$A \rightarrow B(H)$ for some Hilbert space $H$.
In this paper all our operator algebras and spaces will be taken
to be complete.  We shall say that an 
operator algebra is {\em unital} if it has a two-sided
contractive identity.   
In the present paper
we are concerned with operator algebras
with a one-sided (usually right) contractive approximate
identity. 
We shall 
abbreviate `right (resp. left) contractive approximate  identity'
to `r.c.a.i.' (resp. `l.c.a.i.').   
In \S 2 of our paper we 
consider a certain `transference principle',
which can allow one to deduce many general results about 
 operator algebras  with r.c.a.i., from results about
left ideals in a $C^*$-algebra (see the companion 
paper \cite{Bid1}).  Namely there is an important 
left ideal ${\mathfrak J}_e(A)$ of a $C^*$-algebra $\E(A)$, which 
is associated to any 
such operator algebra.  We call ${\mathfrak J}_e(A)$ the `left ideal 
envelope' of $A$.   This is analoguous to 
the case of operator algebras with 2-sided c.a.i., 
which are largely studied these days
in terms of a certain $C^*$-algebra, namely
the $C^*-$envelope.   

    

In \S 3  we study the `left multiplier operator algebra' 
  of an  operator algebra $A$ with r.c.a.i. (which will be 
a symmetrical theory to that of $RM(A)$ for an 
operator algebra $A$ with l.c.a.i).   
 This theory does not work out quite 
as nicely as the case when $A$ has a 
two-sided c.a.i., unless $A$ is a
left ideal of a $C^*$-algebra.  We present some 
examples to show that some of the obvious candidate 
descriptions of this multiplier algebra are defective.      
A main motivation for studying multiplier algebras is because of
their intimate connection with the extremely important notion of the 
`unitization' of a nonunital algebra.  A good `unitization' procedure
should reduce many problems about  nonunital algebras to the well
understood unital case.   Although we are not at this point able
to completely analyze the `unitization' procedure for algebras
with one-sided c.a.i., we at least do some of the groundwork 
here.   


We end the introduction with some
more notation, and some background results which will
be useful in various places.  
We reserve the letters $H, K$ for Hilbert spaces,
and $J$ for a left ideal of a $C^*$-algebra.
If $A$ is an algebra then we write $\lambda
: A \rightarrow Lin(A)$ for the  canonical `left regular
representation' of $A$ on itself.
If $S$ is a 
subalgebra of $A$ then 
the left idealizer of $S$ is the 
subalgebra $\{ x \in A : x S \subset S \}$ of $A$.
Note $S$ is a left ideal in this subalgebra, whence the name.
Similarly for the right idealizer; the (2-sided) 
idealizer is the intersection of the left and right idealizer.

By a `representation' $\pi : A \rightarrow B(H)$ of 
an operator algebra $A$ we shall
mean a completely contractive homomorphism.  If
 $A$ has r.c.a.i. and if we say that 
$\pi$ is {\em nondegenerate}, then at the very least  
we mean that 
$\pi(A) H$ is dense in $H$.  Note that this last 
condition does not imply
in general
 that $\pi(e_\alpha) \zeta \rightarrow \zeta$ for $\zeta \in H$,
where $\{ e_\alpha \}$ is the r.c.a.i.,
as one is used to in the two-sided case.  

\begin{lemma} \label{mano}  Let $A$ be an operator algebra with
either a l.c.a.i. or a r.c.a.i..  Suppose that $\pi :
A \rightarrow B(H)$ is a completely contractive
representation,
set $K$ to be the norm closure of $\pi(A) H$, and let
 $\pi' : A \rightarrow B(K)$ be the induced representation.
Then $\pi'$ is a
completely contractive homomorphism, and
it is nondegenerate in the sense that $\pi'(A) K$ is dense
in $K$.  Moreover if $A$ has a r.c.a.i. and
$\pi$  is isometric (resp. completely isometric, 1-1),
then $\pi'$ is also  isometric (resp. completely isometric,
1-1).
 
Consequently, if $A$ is a subalgebra of $B(H)$ with a
r.c.a.i.,  then we may view $A$ as a nondegenerate
subalgebra of $B(K)$, where $K$ is the closure of $A H$.
\end{lemma}
 
\begin{proof}  That $\pi'$ is a completely contractive
homomorphism is clear.  Clearly $\pi'(A) K \subset K$, however
if $a \in A$ and $\zeta \in H$ then $\pi(a) \zeta =
\lim \pi'(a) \pi(e_{\alpha}) \zeta \in
\pi'(A) K$, if $A$ has a r.c.a.i..   Taking linear combinations
of such expressions $\pi(a) \zeta$, and limits,
shows that any element of $K$ is in the closure of $\pi'(A) K$.
So $\pi'$ is nondegenerate in this case. A similar argument
holds if $A$ has a l.c.a.i..
 
To see the statement here about the isometry note that for $a \in A$
we have
$$\Vert a \Vert  = \Vert \pi(a) \Vert = \sup \{ \Vert \pi(a) \zeta \Vert :
\zeta \in Ball(H) \} =
\sup \lim_\alpha
 \{ \Vert \pi'(a) \pi(e_\alpha) \zeta \Vert :
\zeta \in Ball(H) \} \;.$$
However the right hand side is dominated by
$\Vert \pi'(a) \Vert$, so that $
\pi'$ is an isometry.  A similar argument works for a
complete isometry.  The 1-1 assertion is easier.
\end{proof}
 
\vspace{4 mm}
 
If an operator algebra has only a l.c.a.i.
then the
 `isometric' assertions of the last result are not true
in general.  For a counterexample consider $A = R_n$,
the subalgebra of $M_n$ supported on
the first row.   

We will use without comment several very basic facts from
$C^*$-algebra theory (see e.g. \cite{Ped}), such as
the basic definitions
of the left multiplier algebra $LM(A)$,
and the multiplier algebra $M(A)$, of a $C^*$-algebra.
As a general reference for operator spaces the reader might consult
\cite{ERbook,Pis} or \cite{P}.
We write $\; \hat{ } \; : X \rightarrow X^{**}$ for the 
canonical map, this is a complete isometry if $X$ is an
operator space, and is a homomorphism if $X$ is an operator 
algebra, giving the second 
dual the Arens product \cite{BonsallDuncan}.

We will often consider the basic examples $C_n$ (resp. $R_n$)
of operator algebras with right (resp. left) identity of norm 
1; namely the $n \times n$ matrices
`supported on' the first column (resp. row).  This is a
left (resp. right) ideal of $M_n$, and has the projection
$E_{11}$ as the 1-sided identity. 

If $X$ and $Y$ are subsets of an operator algebra we usually write 
$X Y$ for the {\em norm closure} of the set of finite 
sums of products $x y$ of a term in $X$ and a term in $Y$.
For example, if $J$ is a left ideal of a $C^*$-algebra $A$,
then with this convention $J^* J$ will be
 a norm closed $C^*$-algebra.  This convention extends to 
three sets, thus $J J^* J = J$ for a left ideal of a $C^*$-algebra
as is well known.  These facts follow easily from the 
well known results that
a (norm closed) left ideal $J$ in a $C^*$-algebra has a
positive right contractive approximate
identity.   Also $J \cap J^* = J^* J \subset J \subset J J^*$,
 so that $J$ is a left ideal of $J J^*$.  
 
For the purposes of this paper
we will define a {\em triple system} to be a norm closed
subspace $X$ of a
$C^*$-algebra such that $X X^* X \subset X$.  By `subspace' we 
will allow for example spaces such as $B(K,H)$, regarded 
as the `1-2-corner' of the $C^*$-algebra 
$B(H \oplus K,H \oplus K)$ in the 
usual way.  It is well known that 
triple systems are `the same thing'
as Hilbert $C^*$-modules, although there is a slight difference 
of emphasis in the two theories; 
the important 
structure on a triple system is the `triple product' 
$x y^* z$.  A `triple subspace' is a 
norm closed vector subspace of a triple system which is closed 
under this triple product.
If $X$ is a triple system
then $X X^*$ and $X^* X$ are $C^*$-algebras,
which we will call the left and right $C^*$-algebras of $X$ respectively,
and $X$ is a $(X X^*) - (X^* X)$-bimodule.
  A linear map
$T : X \rightarrow Y$ between triple systems is a {\em triple 
morphism} if $T(x y^* z) = T(x) T(y)^* T(z)$ for all 
$x,y,z \in X$.   Triple systems are operator spaces, and 
triple morphisms are completely contractive, and indeed are 
completely isometric if they are 1-1 (see \cite{Ham2},
this is related to results of Harris and Kaup).
A completely
isometric surjection between TRO's is a triple morphism.  This 
last result might date back to around 1986, to Hamana and 
Ruan's PhD thesis independently.
See \cite{Ham2} or \cite{BShi} A.5 for a proof.

Several times in \S 2 we will refer to 
Hamana's triple envelope
$\T(X)$ of an operator space $X$.  This theory 
may be found in \cite{Ham2}, although we review the
 construction of $\T(X)$ briefly in \S 1 below.   
An alternative account of Hamana's theory of 
the triple envelope is given in \cite{BShi}, 
particularly the introduction and Appendix A there.
The space
$\T(X)$ may be viewed as a triple system, and there is a canonical 
complete isometry $i : X \rightarrow \T(X)$ such that 
there is no nontrivial
triple subspace of $\T(X)$ containing $i(X)$.  We
 write $\E(X)$ and ${\mathcal F}(X)$ for the
left- and right- $C^*$-algebras of $\T(X)$ respectively.     

Next we recall 
the left multiplier algebra
$\M_\ell(X)$ of an operator space $X$.   This is a unital
operator algebra, which may be viewed as a subalgebra of 
$CB(X)$ containing $Id_X$, but with a different (bigger in general)
norm.   Here $CB(X)$ is the space of completely bounded 
linear maps on $X$.
Our first definition of $\M_\ell(X)$ from 
\cite{BShi} was as the linear maps $X \rightarrow X$ corresponding to 
elements $T$ in $LM(\E(X))$ such that $T i(X) \subset i(X)$.
The norm  on $\M_\ell(X)$ is the 
$LM(\E(X))$ norm, and similarly for matrix norms.
There are several other
equivalent definitions of $\M_\ell(X)$ given in
\cite{BShi,BEZ,BPnew} - it is best to view operator space multipliers 
as a `sequence of equivalent definitions'.   In any particular context
one or other of these definitions may be more appropriate.  We 
will make much use the following two simple results from \cite{Bid1}:          
 
\begin{lemma}  \label{lemu}  Let $A$ be an operator algebra
with a  r.c.a.i..  Then the canonical `left regular
representation' of $A$ on itself yields completely contractive
embeddings (i.e. 1-1 homomorphisms)
$$A \hookrightarrow \M_\ell(A) \hookrightarrow CB(A) \; \; ,$$
and the first of these embeddings, and their composition,
are completely isometric.
\end{lemma}
 
\vspace{4 mm}
 
We remark that the canonical
 embedding $\M_\ell(A) \hookrightarrow CB(A)$,
where $A$ is an operator algebra with r.c.a.i.,
 is not in general completely isometric, or even isometric
(see Example \ref{Exmdb}).   This has implications for our theory
of multipliers in \S 3.

\begin{lemma} \label{pid2} 
Suppose that $a \in B(H)$, and
$\{ e_\alpha \}$ is a net of contractions in $B(H) $ such that
$a e_\alpha \rightarrow a$.  Then
$a e_\alpha e_\alpha^* \rightarrow a$,
$a  e_\alpha^* e_\alpha \rightarrow a$, and
$a  e_\alpha^* \rightarrow a$.
\end{lemma}
 
%
 
\section{The ideal envelope}


We begin by sketching Hamana's construction of the injective and 
the triple envelope of operator spaces.  The reader who is not 
familiar with this will need to consult \cite{Ham2} for 
more details and notation; unfortunately
the material below will be a little technical for 
those not versed in `Hamana-theory'.   See also
the introductions of \cite{BPnew} and \cite{BShi}). 
Let $X\subset B(H)$ be 
an operator space, and consider the Paulsen 
operator system
$${\mathcal S}(X):=\left[
\begin{array}{cc}
\Co 1_H & X \\
X^* & \Co 1_H
\end{array}
\right]
\subset M_2(B(H)) \; .$$   
One then takes a minimal (with respect to a certain 
ordering) completely positive 
${\mathcal S}(X)$-projection $\Phi$ on
$M_2(B(H)).$  
By a well-known result of Choi and Effros,
Im$ \;\Phi$  is a $C^*$-algebra with the multiplication $\odot$ defined by
$\xi \odot \eta := \Phi(\xi \eta)$ for $\xi, \eta \in$ Im$\Phi,$ and the
other algebraic operations and norm are the usual ones.   
One may write   
\begin{center}
Im$ \; \Phi = \left[
\begin{array}{cc}
I_{11} & I(X) \\
I(X)^* & I_{22}
\end{array}
\right] \; \subset M_2(B(H)) ;$
\end{center}
and $I(X)$ is an {\em injective envelope}
of $X,$ and $I_{11},$ $I_{22}$ are 
injective C*-algebras.    We sometimes write 
$I_{11}$ as $I_{11}(X)$, say, to emphasize the dependence 
on $X$.  

We continue to think of 
$I_{11},$ $I_{22},$ $I(X)$ and $I(X)^* $
as subsets of $B(H)$, however the operation
$\odot$ induces new products between elements of
$I_{11},$ $I_{22},$ $I(X)$ and $I(X)^*.$   To distinguish these
multiplications from the original product on $M_2(B(H))$ we write
the new products as $\circ$.  

By a well known trick one may also decompose
$$\Phi=\left[
\begin{array}{cc}
\psi_1 & \phi \\
\phi^* & \psi_2
\end{array}
\right].$$ 

Also one may write the $C^*$-subalgebra of $Im \; \Phi$ (with the 
new product) generated by 
\begin{center}
$\left[
\begin{array}{cc}
0 & X \\
0 & 0
\end{array}
\right]      $
\end{center}
as 
\begin{center}
$\left[
\begin{array}{cc}
{\mathcal E}(X) & {\mathcal T}(X) \\
{\mathcal T}(X)^* & {\mathcal F}(X)
\end{array}
\right]
\subset
\left[
\begin{array}{cc}
I_{11} & I(X) \\
I(X)^* & I_{22}
\end{array}
\right] \; .  $
\end{center} 
This defines $\T(X)$, it is clearly a
$C^*$-module or triple system (viewed
as the 1-2-corner of 
the $C^*$-subalgebra just introduced,
 its triple product is $x \circ y^* \circ z$).  
Indeed the span of
expressions of the form $a_1 \circ a_2^* \circ 
a_3 \circ  a_4^* \circ 
\cdots \circ a_{2n+1}$,
for $a_i \in X$, are dense in $\T(X)$.   
Thinking of 
$\T(X)$ as a triple system, and letting $i$ be the canonical 
map $X \rightarrow \T(X)$, we say that $(\T(X),i)$ is
a {\em triple envelope} of $X$.
                                        
 \begin{lemma}
\label{product}
Let A be a subalgebra of $B(H)$, and suppose that A has a r.c.a.i.
$\{ e_{\alpha} \}.$ Then
for all $a, b\in A$ we have (taking $X = A$ and using the notation above)
$$ab=\lim_{\alpha} a\circ
e_{\alpha}^*\circ b \; ,  \; \text{and} $$
$$\psi_1(a) = \lim_\alpha a \circ e_\alpha^*  \in I_{11} \; . $$
\end{lemma}
 
\begin{proof}  Using the notation above, and Lemma \ref{pid2},
we have 
$$\psi_1(a) = \lim_\alpha \psi_1(a e_\alpha^*) = \lim_\alpha  a \circ e_\alpha^*
 \; , $$         
the last step by definition of the $\circ$ product.  Similarly                                          
$$ab=\phi(ab)=\phi(\lim_{\alpha} ae_{\alpha}^*b)
=\lim_{\alpha} a\circ e_{\alpha}^* \circ b . $$ 
\end{proof}  

We recall that an equivalent definition of $\M_\ell(X)$ given in 
\cite{BPnew} was as $IM_\ell(X) =
\{ y \in I_{11} : y \circ X \subset X \}$.
Using the last lemma we see that if $A$ is an operator algebra
with r.c.a.i., then $\psi_1(a) \in IM_\ell(A)$ for all 
$a \in A$.  Thus if $j$ is $\psi_1$ restricted to $A$,
then $j : A \rightarrow IM_\ell(A)$.  By the last lemma,
$j(a) \circ b = ab$ for all $a,b  \in A$, so that $j$ corresponds 
to the left regular representation $\lambda : A
\rightarrow \M_\ell(A) \subset CB(A)$.

\begin{proposition} \label{notfi}  Let $A$ be an operator 
algebra with r.c.a.i..
\begin{itemize}
\item [(1)] Under the same assumptions as  Lemma
\ref{product}, and notation
as above,
$$IM_{\ell}(A)=\{T\in IM_{\ell}(A) :
 T \circ j(A) \subset j(A) \}.$$
\item [(2)]  Elements of
$\M_\ell(A)$, considered as maps on $A$,
are right $A$-module maps.  That is, $\M_\ell(A) \subset CB_A(A)$ as
sets.   Also, $\M_r(A) \subset \; _ACB(A)$ as sets. 
\end{itemize}  
\end{proposition}
 
\begin{proof}
(1) Let $a \in A$, then         
$T\circ j(a)
= \lim_\alpha T\circ a\circ e_{\alpha}^*$,
which equals $j(T \circ a)$
if $T \in IM_{\ell}(A)$.  Conversely,
if $T\circ j(a) = j(a') $ for an $a' \in A$, 
then by the last formula and \ref{product} we
have $T\circ j(a) \circ e_\beta =  T\circ a e_\beta
= a' e_\beta$, so that $T\circ a = a' \in A$.   

(2)  
For the last assertion, if $T \in \M_r(A)$
then there exists $y \in I_{22}$ such that 
$A \circ y \subset A$, and $T$ is the
map $T(a) = a \circ y$.  Hence for $c \in A$,  
$$T(ca) = (ca) \circ y = \lim_\alpha c \circ e_{\alpha}^*\circ
a \circ y =  c T(a) \; \; ,$$
using \ref{product} twice.  The other is similar.
\end{proof}

\vspace{4 mm}

In the statement of the following theorem, and in 
keeping with the usual presentation of Hamana theory
(in this case the triple envelope), we forget the 
original product on $B(H)$ mentioned in the above discussion.
That is, the products in the  statement below are  
the $\circ$ operation in the above discussion.   However 
in the {\em proof} we go back to the $\circ$ notation.  
 
\begin{theorem} \label{Thlem}
Let $A$ be an operator algebra with r.c.a.i. $\{ e_\alpha \}$.
Suppose that $(\mathcal{T}(A),i)$ is a triple envelope 
for $A$, and let $\mathcal{E}(A) = \mathcal{T}(A) 
\mathcal{T}(A)^*$.   The map
$\psi : \mathcal{T}(A)\longrightarrow 
\mathcal{E}(A)$  defined by
$\psi (x):=\lim_{\alpha} x i(e_{\alpha})^*$ is a 
well defined complete
isometry, and $\psi \circ i$ is a homomorphism $j$ on A.   The range of 
$\psi$ is a left ideal $J$ of the $C^*$-algebra
$\mathcal{E} (A)$, and $\psi$ is also a
triple morphism onto $J$.   Thus $(J,j)$ is another triple 
envelope of $A$.
Moreover $J J^* = {\mathcal E} (A),$ and any 
r.c.a.i for $A$ is taken by $j$ to a r.c.a.i. 
for $J$.   Finally, $\psi(x) \psi(y)^*
= x y^*$ for any $x, y \in \mathcal{T}(A)$.
\end{theorem}
 
\begin{proof}  By `abstract nonsense' we can assume that 
$(\mathcal{T}(A),i)$ is the triple envelope considered 
above \ref{product} (in short because
the statement we are attempting
to prove is invariant under the notion we called 
$A$-isomorphism in \cite{BShi} Appendix A). 
We use the notation of \ref{product} and above \ref{product}.
Since elements spanning a dense subset of 
$\mathcal{T}(A)$ are $\circ$ products ending with a term in $A$,
it follows from Lemma \ref{product} that $
\lim_{\alpha} x \circ e_{\alpha}^*$ exists in $\mathcal{E}(A)$
for any $x \in \mathcal{T}(A)$.  Clearly 
$\psi$ is a complete contraction extending the map $j$ 
introduced above \ref{notfi}.
If $x\in A,$ then $$\| \psi (x)\| =\lim_{\alpha} \|
x\circ e_{\alpha}^*\| \geq \lim_{\alpha} \|
x\circ e_{\alpha}^* \circ e_{\beta} \|
= \|x e_{\beta} \|$$
using \ref{product}.  Taking the limit shows that 
$\psi$ is isometric on $A$.  A similar argument using 
associativity of $\circ$ shows that $\psi$ is isometric on $\mathcal{T}(A)$.
Similarly, $\psi$ is a complete isometry.
If $a, b\in A,$ then 
$$\psi(ab)=\lim_{\alpha} (ab)\circ
e_{\alpha}^*
=\lim_{\alpha}(\lim_{\beta}a\circ e_{\beta}^*\circ
b)\circ e_{\alpha}^*
=(\lim_{\beta}a\circ e_{\beta}^*)\circ (\lim_{\alpha}
b \circ e_{\alpha}^*)=\psi (a)\circ \psi (b).$$
 Note also that
$$\psi(a) \circ \psi(b)^* = 
\lim_{\alpha} \lim_{\beta} a\circ e_{\alpha}^* \circ 
e_{\beta} \circ b^* = \lim_{\beta} (ae_{\beta}) \circ b^*
= a \circ b^* \; \; .$$
Thus in the language of the statement of the 
theorem $\psi(i(a)) \psi(i(b))^* = i(a) i(b)^*$.   
 
By looking at the natural dense subsets of 
$\mathcal{T}(A)$ and $\mathcal{E}(A)$ it is easy 
to argue that the range of
$\psi$ is a left ideal $J$.  By a fact mentioned in the 
introduction this implies that $\psi$ is a triple 
morphism.   Since $j$ is the restriction of 
$\psi$ to $A$, the remaining 
assertions are easy.   One needs to use the fact that 
$\psi$ is a triple morphism, and that $\mathcal{T}(A)$
and therefore consequently $J$, have  dense subsets spanned by
terms which are alternating products as mentioned 
above \ref{product}. \end{proof}
 
\vspace{5mm} 
 
\begin{theorem} \label{Thlemi}  Let $A$ be an operator algebra
with a r.c.a.i., let $(I(A),i)$ be its injective envelope 
 discussed above \ref{product} say,
and let $B$ be the injective unital  
$C^*$-algebra $I_{11}$ (also discussed above).    Then there is an
orthogonal projection $e \in B$, such that $I'(A) = (B e,j)$ is also an
injective envelope for $A$, where $j$ is as above.
For this new injective envelope, $I'_{22} \cong eBe$ and 
$I'_{11} \cong B$.  The completely isometric embedding
$j : A \rightarrow I'(A) \subset B$ is a
homomorphism.  Moreover if $A$ has a right identity $f$ then 
$j(f) = e$.
\end{theorem}  

Before we prove this, we remark that we may define the {\em ideal 
injective envelope} 
of an operator algebra $A$ with r.c.a.i.,
to be $I'(A)$ as in the theorem.  
 
\vspace{4 mm}
 
\begin{proof}   
%
Define $\psi : I_{11} \circ A
\rightarrow I_{11}$ by $\psi(x) = \lim_\alpha
x \circ e_\alpha^*$.  By \ref{product} this limit exists and 
$\psi$ is clearly completely contractive.  Note that 
$\psi$ is also a left $I_{11}$-module map, so that by a
result in \cite{BPnew} we may extend $\psi$ to a
completely contractive left $I_{11}$-module map 
 $\hat{\psi} : I(A) \rightarrow I_{11}$.   The restriction 
of this map to $A$ was the map called $j$ above.
By Hamana's
`essential property' for $I(A)$ \cite{Ham2}, $\hat{\psi}$ is a 
complete isometry.  By \cite{BPnew} 2.7, there is
a $v \in Ball(I(A))$ such that $\hat{\psi}(x) = x \circ v^*$,
for all $x \in I(A)$.  Also its range is a left ideal of 
$I_{11}$, so that $\hat{\psi}$ is a triple morphism 
by a fact at the end of the introduction.

Next define a map 
$I_{11} \circ j(A) \subset I_{11} \rightarrow I(A)$ by
$z \mapsto \lim_\alpha z \circ e_\alpha$.   To see that 
this limit exists note that for $x \in I_{11}$ and $a \in A$
we have by the type of calculations found in 
the proof of \ref{Thlem}, that $\lim_\alpha x \circ j(a) \circ e_\alpha
= x \circ a$.  So the limit exists on $I_{11} \circ j(A)$, and moreover
this map is exactly the map $x \circ j(a) \mapsto 
x \circ a$.   It is thus a completely contractive left 
$I_{11}$-module map, and extends by \cite{BPnew} to a
completely contractive left $I_{11}$-module map
$\mu : I_{11} \rightarrow I(A)$.  There clearly exists
$w \in Ball(I(A))$ such that $\mu = - \circ w$.  Now notice 
that the composition $\mu \hat{\psi}$ is a completely 
contractive map $I(A) \rightarrow I(A)$, which restricts to the 
identity map on $A$.   By Hamana's rigidity property \cite{Ham2},
$\mu \hat{\psi} = Id$ on $I(A)$.  Thus $\hat{\psi}$ is a
complete isometry, $\mu$ is onto.  Moreover since 
$x \circ (v^* \circ w) = x$ for all $x \in I(X)$, it follows 
from 
 \cite{BPnew} Corollary 1.3, that $v^* \circ w$ is the 
identity of $I_{22}$.  Thus by a well know operator theory fact,
$v = w$.    We define $e = v \circ v^* \in I_{11} = B$, this 
is an orthogonal projection.   Moreover $Ran \; \hat{\psi} =
B e$, as may be seen easily from the above.  It follows 
immediately that  $(B e,j)$ is an 
injective envelope for $A$, where $j(a) = a \circ v^*$ (which coincides
with the map $j$ of previous results), which is a 
homomorphism.  

Note that by \cite{BPnew} Theorem 1.8, $I'_{11} \cong M(I'(X) I'(X)^*)
= M(I(X) I(X)^*) \cong I_{11}$, and $I'_{22} \cong M(eBe) = eBe$. 

Finally, if $A$ has right identity $f$ then we may take 
$v$ above to be $f$, so that $j(f) = f \circ f = e$.   
\end{proof} 
 
\vspace{4 mm}

It is clear from \ref{Thlem}
that for an operator algebra $A$ with r.c.a.i., $\T(A)$ 
may be taken to be a left ideal in a $C^*$-algebra, and this
left ideal therefore  
possesses  Hamana's universal property of the triple 
envelope (see \cite{Ham2}
or \cite{BShi} Appendix A).  However examining this property in 
this case gives an interesting refinement, by viewing 
$(\T(A),j)$ as an {\em left ideal extension} of $A$.   

To be 
more specific we say that a pair $(J,i)$ consisting 
of a left ideal $J$ in a $C^*$-algebra, and a completely isometric
homomorphism $i : A \rightarrow J$, is 
a {\em left ideal extension} of $A$
if $i(A)$ `generates $J$ as a triple system'.  That is, the span of 
expressions of the form $i(a_1) i(a_2)^* i(a_3)  i(a_4)^* 
\cdots i(a_{2n+1})$,
for $a_i \in A$, are dense in $J$.   It follows from this
that $\{ i(e_\alpha) \}$ is a r.c.a.i. for $J$ if 
$\{ e_\alpha \}$ is a r.c.a.i. for $A$.  Then it is clear from 
\ref{Thlem} that $(\T(A),j)$ may be taken to 
be a left ideal extension of $A$, 
and it is the `minimal such'.   The new point in the theorem below
is that $\tau$ may be chosen to be a homomorphism:

\begin{corollary} \label{upie}  Let $A$ be an 
operator algebra with r.c.a.i..  Then there 
exists a left ideal extension $({\mathfrak J}_e(A),j)$ of $A$, 
with ${\mathfrak J}_e(A)$ a left ideal in a $C^*-$algebra
$\E(A) = {\mathfrak J}_e(A) {\mathfrak J}_e(A)^*$, such that for 
any other left ideal extension $(J,i)$ of $A$, there exists 
a (necessarily unique and surjective)
completely contractive homomorphism $\tau: J 
\rightarrow {\mathfrak J}_e(A)$, which is also a
triple morphism, such that $\tau \circ i = j$.  Thus 
${\mathfrak J}_e(A)/(ker \; \tau)  \cong J$ completely isometrically
homomorphically (i.e as operator algebras) too.
Moreover $({\mathfrak J}_e(A),j)$ is unique in the following sense:
given any other $(J',j')$ with this universal property,
then there exists a surjective completely isometric 
homomorphism $\theta : {\mathfrak J}_e(A) \rightarrow J'$ such that 
$\theta \circ j = j'$.
\end{corollary}

\begin{proof}  The uniqueness is fairly obvious and 
standard.  The existence follows from 
\ref{Thlem} (setting ${\mathfrak J}_e(A) = (J,j)$ there),
together with Hamana's universal property 
for the `triple envelope' mentioned above. Indeed we have:
 $$\tau(x y) = \lim \tau(x i(e_\alpha)^* y)
= \lim \tau(x) j(e_\alpha)^* \tau(y) = \tau(x) \tau(y)$$
for $x, y \in J$, using    
Lemma \ref{pid2} and the observations above.
\end{proof}

\vspace{4 mm}

We call $({\mathfrak J}_e(A),j)$ the {\em left ideal envelope} of 
$A$, and continue to write $\E(A) = {\mathfrak J}_e(A) {\mathfrak J}_e(A)^*$.
Again, $j$ will be called the {\em Shilov embedding
homomorphism}.

%
%
 

\begin{corollary} \label{usf}  Let $A$ be an operator algebra
with r.c.a.i., and let $j : A \rightarrow \E(A)$ be the 
Shilov embedding homomorphism.
Then $\M_\ell(A) \cong \{ T \in LM(\E(A)) : T j(A) \subset j(A) \}$
completely isometrically isomorphically.  \end{corollary} 
 
\vspace{4 mm}
      
This last corollary above
is a restatement of \ref{notfi}, using
the definition of $\M_\ell(A)$ given in the introduction.
 
\section{The left multiplier algebra of an algebra with right 
approximate identity}

The readers first thought here might be that this case must be
analogous to
the case of the left multiplier algebra
of an operator algebra with 
l.c.a.i., which was treated in \cite{Bid1}.
However in fact the theory is rather different.
Part of this may be seen by considering the basic example
of $A = C_n$ (the $n \times n$ matrices
`supported on' the first column).  
It is fairly clear that
its left multiplier operator algebra `should be'
$M_n$.

After a little thought about operator algebras 
with r.c.a.i., one comes up with the following list
(which we shall refer to throughout this section) of
 `possible candidate definitions' for $LM(A)$ in this
case (all of which contain $A$ completely 
isometrically):   
\begin{enumerate}
\item [(I)] $\{ x \in A^{**} : x \hat{A} \subset \hat{A} \}$,
\item [(II)] $\{ T \in B(H) : T \pi(A) \subset \pi(A) \}$, where
$\pi : A \rightarrow B(H)$ is a  
completely isometric nondegenerate representation,
\item [(III)] $\M_\ell(A)$, 
\item [(IV)] $CB_A(A) \; \; \; \; $ (right module maps), 
\item [(V)] the `maximum essential left multiplier extension of $A$'.
\end{enumerate}
We will spell out later what is meant in this setting by 
 `essential left multiplier extension'.
One can rule out (I) fairly quickly as a plausible candidate
since it is not unital (e.g. take $A = C_n$).
However one would hope that most of the other four items are
completely isometrically isomorphic.   Unfortunately, most of the
equivalences amongst these items 
break down if $A$ is a general operator algebra
with r.c.a.i..
Nonetheless one of these five candidates will
emerge from our study below as the `winner', namely as the
appropriate version of
$LM(A)$ in this case.  This seems quite interesting.
              
In the special case of a  
left ideal $J$ in a $C^*$-algebra things are much better.
Most of these five items  are then
 completely isometrically isomorphic.    Indeed
in this case we shall see that most of these coincide with the 
usual left multiplier algebra
$LM(J J^*)$ of the $C^*$-algebra $J J^*$.  
One may view this as being predicted by the theory
of Hilbert $C^*$-modules.  See also \cite{KP}.  Indeed 
the analogous `left multiplier $C^*$-algebra' $LM^*(J)$,
which may be defined to be $M(J J^*)$,
or equivalently $A_l(J)$ in the notation of \cite{BShi},
also has a very satisfactory theory which has 
essentially been done in \cite{KP}.  Hence we 
shall not mention $LM^*(J)$ again here, leaving any details
of its theory to the interested reader.
  
In this special case of a left ideal $J$ in a $C^*$-algebra
we will take  the term 
`essential left multiplier extension'  used 
in (V), to be defined as in
the proof of (4) in the main theorem 
in \cite{Bid1} \S 4.  Of particular 
interest, in (II),
are the {\em ideal representations} of $J$
discussed in \cite{Bid1}, that is,
the completely isometric representations $\pi : J
\rightarrow B(H)$ which 
are restrictions of faithful 
*-representations $J J^* \rightarrow B(H)$ (or equivalently,
the completely isometric representations of $J$ which are also 
triple morphisms). 
 


\begin{theorem} \label{lmli}  Let $J$ be a left ideal of a 
$C^*$-algebra.   The algebras in (III)-(V) are all
completely isometrically isomorphic to
$LM(J J^*)$.   Any algebra in (II)-(V) contain
 the algebra in (I) completely isometrically isomorphically
as a (proper, in general) subalgebra.  If, further,
 in (II) we only consider 1-1 nondegenerate ideal representations 
$\pi$, 
then the algebras in (II)-(V) are all
completely isometrically isomorphic.
 Also $B_J(J) = CB_J(J)$ isometrically. 
\end{theorem}

\begin{proof}  In this 
case $J$ is a $C^*$-module, and $\M_\ell(J)$ may easily be seen, from
\cite{BShi} A.4 if necessary, to be
$LM(J J^*)$.   Clearly $CB_J(J) \subset CB_{J^*J}(J)$, and conversely,
if $T \in CB_{J^*J}(J)$ then for $x, y \in J$ we have
using Lemma \ref{pid2}
that 
$T(xy) = \lim T(x e^*_\alpha y) = \lim T(x) e^*_\alpha y
= T(x) y$, if $\{ e_\alpha \}$ is the r.c.a.i. in 
$J$.   Thus $CB_J(J) = CB_{J^*J}(J)$.
By a theorem of Lin (\cite{Lin} 1.4)
$B_{J^*J}(J) =  LM(J J^*)$, and the operator space version of 
this is true too (see \cite{Bna}), that is
 $CB_{J^*J}(J) = LM(J J^*)$ completely isometrically.
Thus (III) = (IV), and we obtain 
the last assertion of the statement of the theorem too.   
Looking at the definition 
of a
`essential left multiplier extension', and the associated `ordering', below,
or defining these as we did in 
 the proof of (4) in the main theorem
in \cite{Bid1} \S 4, it 
is now clear that (IV) is an essential left multiplier extension,
so that it is the  
biggest such.   That is, (IV) = (V).   

We will defer until \ref{lmrci} the proof that an algebra in (II)
contains the algebra in (I) completely isometrically.
Finally, 
given a nondegenerate
faithful *-representation $\theta : J J^* \rightarrow B(H)$, 
we know from the multiplier theory for $C^*$-algebras (see
e.g. \cite{Ped} 3.1.12)
that 
$LM(J J^*) = \{ T \in B(H) : T \theta(J J^*) 
\subset \theta(J J^*)$.   However since for any ideal 
$J$ we have $J = J J^* J$, if  
$T \in B(H)$, then  $T \theta(J) \subset \theta(J)$
 if and only if  $T \theta(J J^*) \subset \theta(J J^*)$.
Thus the algebra in (II) equals $LM(J J^*)$ too in this
case.
\end{proof} 

\vspace{4 mm}
 
We have now completed our discussion of $LM(J)$ for a left ideal 
$J$ in a $C^*$-algebra, and turn to the 
more general case of an arbitrary operator algebra 
$A$ with r.c.a.i..   Since this 
case is the most complicated, we will be a little more formal
and rigorous in our presentation below.   

We begin by defining a {\em  left multiplier
extension} of $A$
to be a pair $(B,\pi)$ consisting of an operator algebra
$B$ with an identity of norm 1, and a completely isometric
homomorphism $\pi : A \rightarrow B$ such that $\pi(A)$ is a
left ideal of $B$.   We say that $(B,\pi)$
is an {\em essential left multiplier extension} of $A$
if in addition the canonical map
$B \rightarrow B(A)$ is 1-1.     Note that this canonical map
is then a completely contractive homomorphism into
$CB_A(A)$ (viewing $A$ as a right $A$-module).  In the
2-sided case $CB_A(A)$ is itself a left multiplier
extension of $A$, and is therefore the largest
left multiplier extension, but this is not true in general for us
(see example \ref{Exma}).   We next define an ordering
on essential left multiplier
extensions of $A$: namely $(B,\pi) \leq (B',\pi')$ if there
exists a (necessarily unique and one-to-one)
completely contractive homomorphism $\theta :
B \rightarrow B'$ such that $\theta \circ \pi = \pi'$.
We remark that in the theory in  \cite{Bid1} \S 4 or for
left ideals in $C^*$-algebras one may insist that
$\theta$ is completely isometric, but we are not certain
if we can make this requirement in general.  Note that this ordering
says that $B$ may be viewed as a subalgebra of $B'$, but
with a possibly bigger norm.
 
 
We will use the fact that
an essential left multiplier extension $(B,\pi)$
has the following rigidity property: if
$\theta : B \rightarrow B$ is a completely contractive
homomorphism extending the identity mapping
on $\pi(A)$, then $\theta$ is the identity mapping
(To see this note that
$\theta(b) \pi(a) = \theta(b\pi(a)) = b\pi(a)$,
so that $\theta(b) - b$ is in the kernel of the canonical
map $B \rightarrow CB(A)$).                                

We say that two left multiplier extensions
$(B,\pi)$ and $(B',\pi')$ are
{\em $A$-equivalent} if there
exists a unital completely isometric isomorphism $\theta :
B \rightarrow B'$  with $\theta \circ \pi = \pi'$.  Note that
this is an equivalence relation,
and that  `$\leq$' induces a well defined ordering on the
equivalence classes.   It follows that if
there exists  a maximum essential
 left multiplier extension of $A$, then it is
unique up to $A$-equivalence.
We remark further that if
two left multiplier extensions are $A$-equivalent,
and if one is essential, then so is the other.  We leave
the proof of this as an easy exercise.       

One more definition: we say that a representation $\pi : A
\rightarrow B(H)$ is a completely 
isometric {\em Shilov representation}, if 
$\pi$ is the restriction of a 1-1
*-homomorphism $\E(A) \rightarrow B(H)$
(viewing $A \subset \E(A)$ as in 
\ref{upie} say
via the completely isometric homomorphism 
$j$ there).   Note that the `ideal representations'
 $\pi$ considered 
in Theorem \ref{lmli} happen also to be Shilov 
representations.
 
\begin{theorem} \label{lmrci}  Let  $A$ be an operator algebra
with r.c.a.i..   Then $\M_\ell(A)$ is the (unique up to
$A$-equivalence) maximum 
essential left multiplier extension of $A$.
That is, the algebras in items (III) and (V) 
in the list at the beginning of \S 5 are completely 
isometrically isomorphic.  
Also $\M_\ell(A)$  
contains, and dominates in the ordering above, the 
essential left multiplier extensions described in (II),
which in turn contain the algebra in 
(I) completely isometrically isomorphically.
However if $\pi$ in (II) is a Shilov representation of $A$
then 
the algebra in (II) is completely isometrically isomorphic to
(III) and (V).
The algebra in (IV), namely $CB_A(A)$,
is not in general an operator algebra, and its subalgebra
which corresponds to 
the algebra in (III) (and (V)) does not correspond isometrically.
That is, the canonical 1-1 homomorphic embedding of $\M_\ell(A)$ into 
$CB_A(A)$ (or even $B_A(A)$) need not be 
isometric.
\end{theorem}

\begin{proof}   
That the canonical inclusion of (III) in (IV) is not isometric
may be seen in Example \ref{Exmdb}.   That $CB_A(A)$ is 
not in general an operator algebra may be seen in Example \ref{Exma}.

We have seen in \ref{notfi} that
$(\M_\ell(A),\lambda)$
is an essential left multiplier extension.
It is evident from the fundamental
properties of $\M_\ell(A)$ (from 
Theorem 4.10 (1) of \cite{BShi} say)
that given any essential
left multiplier extension $(B,\pi)$ of $A$, there is a
canonical completely contractive homomorphism
$\sigma : B \rightarrow \M_\ell(A)$.  It is obvious  that
via this homomorphism $(B,\pi) \leq (\M_\ell(A),\lambda)$. 
So $(\M_\ell(A),\lambda)$ is the maximum essential
left multiplier extension.  

If  $\pi : A \rightarrow B(H)$ is a
completely isometric
`nondegenerate representation' (by which we mean at least
that $\pi(A) H$ is dense in $H$), then
$\{ T \in B(H) : T \pi(A) \subset \pi(A) \}$ is fairly 
evidently an essential
left multiplier extension of $A$, and hence is dominated
by $\M_\ell(A)$ (by the fact in the previous paragraph).  


We now show any algebra in
(II) contains (I) completely isometrically.  
If $\pi : A \rightarrow B(H)$ is any (not necessarily 
nondegenerate) completely isometric homomorphism,  
consider the following sequence of
completely contractive homomorphisms:
$$A \hookrightarrow A^{**} \overset{\pi^{**}}{\rightarrow}
B(H)^{**} \rightarrow B(H)$$
where the first and last maps here are the canonical
inclusion and projection.   Let $\sigma$ be the composition of the
last two maps, restricted to $\{ x \in A^{**} : x \hat{A}
\subset \hat{A} \}$.  Then $\sigma(\hat{a}) = \pi(a)$ for
$a \in A$, and so for $x$ as above and $a \in A$ we have
$$
\sigma(x) \pi(a) = \sigma(x \hat{a})  \in \sigma(\hat{A}) 
\subset \pi(A)
 \; . $$   Thus $\sigma$ maps into the left idealizer of $\pi(A)$ 
in $B(H)$.  
To see that  $\sigma$ is isometric note that 
$$\Vert \sigma(x) \Vert \geq \Vert \sigma(x) \pi(e_\alpha) \Vert 
= \Vert \sigma(x \hat{e_\alpha}) \Vert = 
\Vert \pi(x \hat{e_\alpha}) \Vert = \Vert x \hat{e_\alpha} \Vert \; .$$
Now $1$ is a weak* limit point of the $\{ \hat{e_\alpha} \}$, and using
the seperate weak* continuity of the product on 
$A^{**}$,  we see  that $\Vert \sigma(x) \Vert \geq \Vert x \Vert$.
Similarly  $\sigma$ is completely isometric. 

Finally we discuss the equivalence of (III) with 
(II) for Shilov representations $\pi$.
Take $\pi : \E(A)
\rightarrow B(H)$ to be any nondegenerate 1-1
*-homomorphism, then it
is easy to see that
$\pi_{|_A}$ is a nondegenerate completely isometric
homomorphism.   By the basic theory of multipliers of
$C^*$-algebras we may view
$LM(\E(A)) = \{ T \in B(H) :  T \pi(\E(A))
\subset \pi(\E(A)) \}$.   From 
\ref{usf} we thus have  $\M_\ell(A)$ completely 
isometrically isomorphic to
 $$\{ T \in B(H) :  T \pi(\E(A))
\subset \pi(\E(A)), \; \text{and} \;
T \pi(A) \subset \pi(A) \} = 
\{  T \in B(H) :  T \pi(A) \subset \pi(A) \}$$
since $j(A) j(A)^*$ generates $\E(A)$ as a $C^*$-algebra.
\end{proof}

\vspace{4 mm}



\begin{proposition} \label{ifoP}
Let $A$ be an operator algebra with r.c.a.i..
Then $CB_A(A)$ is an operator algebra if and only if
$CB_A(A)$ is completely isometrically isomorphic to $\M_\ell(A)$.
\end{proposition}

\begin{proof}  The one direction is trivial since 
 $\M_\ell(A)$ is an operator algebra.
To see the other note that by \ref{lemu} and \ref{notfi}
we have
 $\lambda(A) \subset CB_A(A)$, and $\M_\ell(A)
\subset CB_A(A)$ although the $\M_\ell(A)$-norm 
is possibly larger than the $CB_A(A)$-norm.   
Thus if $CB_A(A)$ is an
operator algebra then $(CB_A(A),\lambda)$ is an essential
left multiplier extension of $A$
a 1-1 completely contractive
homomorphism $CB_A(A) \rightarrow \M_\ell(A)$ (by the
equivalence of (II) and (V) in \ref{lmrci}).
It is easy to see that this homomorphism restricts to
the identity map on $\M_\ell(A)$, whose inverse is
completely contractive as a map $\M_\ell(A)
\rightarrow CB_A(A)$.
\end{proof}

\vspace{3 mm}

Next we consider some examples:

\begin{example} \label{Exma} \end{example}
Let $A \subset M_3$ be the subalgebra supported on the 1-2, 1-3, 2-2, 
and 3-3 entries only:
$$ A \; = \; \left[ \begin{array}{ccl} 0 & * & * \\
0 & * & 0 \\ 0 & 0 & * \end{array} \right] \; \; .$$
 Then $A$ is an
operator algebra with right identity of norm 1.
One may easily compute ${\mathfrak J}_e(A) = \T(A) = I(A) =
M_{3,2}$, and $\E(A) = M_3$, from which it follows that $\M_\ell(A)$ is the 
subalgebra of $M_3$ spanned by $A$ and $I_3$.  Thus $\M_\ell(A)$ is 
5 dimensional.  Indeed in this example, $\M_\ell(A)$ coincides with  what 
we called the `multiplier unitization' of $A$.
On the other hand the algebra given by item (I) in 
the list at the beginning of \S 5, is $A$ 
again, which is 4 dimensional.  A tedious but completely elementary 
algebraic computation shows that $CB_A(A)$, item (IV) on the aforementioned
list, is 8 dimensional.  This shows by \ref{ifoP}  that 
$CB_A(A)$ cannot be completely isometrically isomorphic to an 
operator algebra (in contrast to the 2-sided c.a.i. case). 

\begin{example} \label{Exmdb} \end{example}   It is perhaps true in the 
previous example that $\M_\ell(A)$ is not contained isometrically in 
$CB(A)$, but this seems difficult to check.  Instead we generalize 
this example to an interesting example of an
operator algebra $A$ with right identity of norm 1,
for which $\M_\ell(A)$ is fairly 
clearly  not contained isometrically in
$CB(A)$ (or $B(A)$).  To do this we take a $C^*$-algebra
$B$ with no nontrivial 2-sided
ideals, and a closed subspace $X \subset B$
which generates $B$ as a triple system.  
 To be more concrete 
one could take $B = M_3$ and 
$X \subset M_3$ as in Example 4.4 of \cite{BShi}.  
Consider $A \subset M_3(B)$ the 
subalgebra supported on the 1-2, 1-3, 2-2,
and 3-3 entries only (i.e. zero elsewhere), which has scalar multiples 
of the identity operator in the 2-2 and 3-3 entries, and has 
operators from $X$ in the 1-2 and 1-3 entries.  If one computes the 
triple system generated by
 $A$ in $M_3(B)$, one sees that one obtains 
$M_{3,2}(B)$ as in the previous example.  As in that example one 
sees that $\E(A) = M_3(B)$ and $\M_\ell(A)$ is the subalgebra of $M_3(B)$
spanned by $A$ and matrices supported only in the 1-1 entry, and that 
entry may be anything in $\M_\ell(X)$.  Now consider one of these 
last matrices, $T$ say, and let $a$ be its one nonzero 
entry (in $\M_\ell(X)$).  
Viewing $T \in CB(A)$ as a left multiplication map, 
it is clear that its norm or `cb-norm' is the same as the 
norm or `cb-norm' of its 1-1 entry, as a left multiplication map on $X$.
But in Example 4.4 of \cite{BShi} we showed that neither of these numbers
equals $\Vert a \Vert = \Vert T \Vert$ in general.

\vspace{4 mm}



{\bf Remark.}   If $A$ is an operator algebra
with r.c.a.i., then there is a most important element of $A^{**}$, namely
the right identity $E$ there.  A very natural question arises as to whether
$E \in LM(A)$, and more particularly whether $E \hat{A} \subset \hat{A}$?
In fact this is not true in general.  To see this suppose that 
we have an operator algebra $A$ with r.c.a.i. but no l.c.a.i.,
which has property ($\R$) discussed 
in \cite{Bid1}.   If $
E \hat{A} \subset \hat{A}$ then in the notation of  
\ref{notfi} there is a corresponding element of 
$e \in IM_l(A)$ such that $e \circ a = E a$ for $a \in A$.  In particular
for $b$ in the algebra $\R(A)$ defined in
\cite{Bid1} \S 2, we have $e \circ b = b$.    If one chooses  
an $A$ with $I(\R(A)) = I(A)$  (such examples exist)
then by the `rigidity property' of 
the injective envelope \cite{Ham2}, it follows that $e \circ a
= a$ for all $a \in I(A)$.  Thus in particular, for $a \in A$ we have
$E a = a$, from which it follows by \cite{BonsallDuncan} 28.7 that 
$a$ has a l.c.a.i., which is a contradiction.

\vspace{4 mm}
      
{\bf Remark.}  (Unitizations.)   Notice that any time 
one has a (essential) left multiplier extension $(B,i)$ 
of an operator algebra $A$ with r.c.a.i., we get a
unitization of $A$ defined as $i(A) + \Co 1_B$.  
Conversely, any unitization of $A$ should be a
left multiplier extension of $A$.  

Taking 
the multiplier extension $\M_\ell(A)$ here, let us 
call the corresponding unitization
the `Shilov unitization' $A^1$
of $A$; it will be the unitization with the 
smallest norm. 
Of course another 
universal essential left multiplier extension may be 
constructed by considering a left idealizer in the 
direct sum of all essential left multiplier extensions of 
$A$.  We obtain an associated unitization of $A$, this one 
with the largest norm.  Indeed this may also described 
by simply assigning $A^1$ the supremum of all operator algebra
structures on $A^1$ coinciding with the usual one on $A$.

We do not know
in general if there is a simple 
formula for the norm on the unitization in terms of the 
norm on $A$, as there is in the case of an algebra
with two-sided c.a.i..   Of course if the 
span if $Id_X$ and $\lambda(A)$ inside $CB_A(A)$ is
an operator algebra, then this coincides  with 
the `Shilov unitization' mentioned above, by an argument
similar to \ref{ifoP}.
These matters deserve investigation.

\vspace{5 mm}

Finally, we end with some remarks on the 2-sided 
multiplier algebra
$M(A)$ of an operator algebra with l.c.a.i. say.   One natural 
candidate for $M(A)$ might be constructed as follows: take 
any candidate $(D,\mu)$ for $LM(A)$ (as in \cite{Bid1} \S 4), and define 
$M(A) = \{ x \in D : \mu(A) x \subset \mu(A) \}$.  It is easy to 
see that if we take a $(D',\mu')$ which is 
$A$-isomorphic to $(D,\mu)$, then one obtains the same $M(A)$,
up to completely isometric isomorphism, or up to a suitable 
notion of $A$-equivalence).  
Perhaps this
coincides with $M(\Li(A))$ if $A$ has property ($\Li$).  

One would hope that one might get the same 
algebra $M(A)$ by looking at $\{ x \in M_r(A) : x \rho(A) 
\subset \rho(A) \}$ where $\rho$ is the `right regular antirepresentation'
of $A$ on itself, but this is not correct
(consider $R_2$).  In fact this algebra is related to a `big' multiplier
algebra $BM(A)$ of $A$, 
defined as the `double idealizer'
of $j(A)$ in $M(B)$, where $B$ is the injective $C^*$-algebra
in \ref{Thlemi}.                                 
 Another approach to $M(A)$ is via
`double centralizers' as in \cite{PuR} and \cite{Pal} \S 1.2.
There must be some relations 
between all these candidates for $M(A)$, and some of the algebras
discussed at the end of  \cite{Bid1} \S 4 will also play a role.   The 
second author is currently investigating these matters
 to see if there are any satisfactory results here.  
 
\vspace{5 mm}

Acknowledgement:  Several results in this paper  will be included as   
part of the second author's Ph. D. thesis, done
under the direction of the first author during the '00-'01 academic 
year.   He expresses his gratitude to the first author for his kind 
instruction and advising.

\end{document}